

\magnification=\magstep1
\baselineskip=13pt
\def \Box {\hbox{}\nobreak \vrule width 1.6mm height 1.6mm
depth 0mm  \par \goodbreak \smallskip}
\newcount\tno     
\tno=1

\noindent{\hfill November 15, 1995}
\vskip .3cm
\centerline{\bf Evolutions, Symbolic Squares, and Fitting Ideals}
\vskip .2cm
\centerline{by}
\vskip .2cm
\centerline{{\bf David Eisenbud} and {\bf Barry Mazur}}
\vskip .5cm
\hrule
\midinsert
\centerline{\bf Abstract}
\noindent Given a reduced local
algebra $T$ over a suitable ring or field $k$ we
study the question of whether there are nontrivial
algebra surjections
$R\to T$ which induce isomorphisms
$\Omega_{R/k}\otimes T \to \Omega_{T/k}$. Such maps,
called evolutions, arise
naturally in the study of Hecke algebras, as
they implicitly do
in the recent
work of Wiles, Taylor-Wiles, and Flach.  We show that the
existence of non-trivial evolutions of an algebra $T$
can be characterized in terms of the
symbolic square of an ideal defining $T$.  We give a
characterization of the symbolic square in terms of Fitting ideals.
Using this and other techniques we show that certain classes of
reduced algebras
--- codimension 2 Cohen-Macaulay, Codimension 3 Gorenstein,
licci algebras in general, and some others --- admit no nontrivial
evolutions.  On the other hand we give examples
showing that non-trivial evolutions of reduced Cohen-Macaulay
algebras of codimension 3 do exist
in every positive characteristic.
\endinsert
\hrule
\vskip .5cm

{\bf Introduction}
\vskip .2cm

\vskip .1cm
\noindent{\sl Definition.}
Let $\Lambda$ be a ring and let $T$ be
a local
$\Lambda$-algebra essentially of finite type,
that is, a localization of a finitely generated $\Lambda$ algebra.  An
{\it evolution} of $T$ over $\Lambda$
is a local $\Lambda$-algebra $R$
essentially of finite type and a surjection $R\to T$
of $\Lambda$-algebras
inducing an isomorphism
$\Omega_{R/\Lambda_{R}}\otimes T \to \Omega_{T/\Lambda}$.
The evolution is {\it trivial} if $R\to T$ is an isomorphism.
\vskip .1cm

This notion was first formulated (in slightly different language)
by Scheja and Storch [1970] and B\"oger [1971]. In their terminology
a pair of ideals $J \subseteq I$ in a polynomial ring P define an evolution
$P/J \to P/I$ if I is {\it differentially dependent} on $J$, and $P/I$
admits
no nontrivial evolutions (that are homomorphic images of $P$) if $I$ is
{\it differentially basic}.  The idea was studied in a number-theoretic
context in Mazur [1994].

For example if
$f \in P := \Lambda[x_1,\dots,x_n]_{(x_1,\dots,x_n)}$
and $I$ is the ideal
generated by the derivatives $f$, then
$P/I$ is an evolution of $P/(I,f)$,
nontrivial when $f$ is not contained in $I$
(in characteristic zero
this means roughly that $f$ is not quasihomogeneous --- see
Saito [1971]).
This may lead one to think that non-trivial evolutions are everywhere.
On the
other hand, if $T$ is a reduced complete intersection,
generically separable
over $\Lambda$, then every evolution of $T$ is trivial
(this well-known result is generalized below). Moreover, we have been
unable to find any nontrivial evolution
of any reduced algebra in equi-characteristic zero, or of any reduced
algebra which is flat over a discrete valuation ring of mixed characteristic.

In this paper we study the general question of existence of nontrivial
evolutions. First we exhibit an elementary characterization,
found by H. Lenstra,
of those algebras that admit a nontrivial
evolution in terms of their cotangent sequences.
We use this
to connect the existence of evolutions with the following question:
given a prime ideal (or more generally an unmixed ideal $I$ in
a local ring $(P,{\cal M})$, when is it true that the symbolic
square $I^{(2)}$ is contained in the product $I{\cal M}$?  (Definitions
are given below).  We characterize symbolic squares in terms of
the behavior of certain Fitting ideals, and use this characterization
to show that  $I^{(2)} \subseteq I{\cal M}$ in the case of
perfect ideals of codimension 2.  We generalize this relation,
using a result of Buchweitz, to the case of licci ideals.
We then show that the same relation (and in some cases
an even stronger one) holds for several other classes of
ideals.  We exhibit a result of E. Kunz showing that this is
the case for unmixed almost complete intersections.  We show that it
also holds for unmixed quasihomogeneous ideals in characteristic 0.
On the negative side, we give examples of unmixed
quasihomogeneous prime ideals in every positive characteristic
with non-trivial evolutions; these were found through an analysis of
an example by E. Kunz, with the help of Sorin Popescu.

We may pose a question which has a somewhat elementary appearance,
but for  which an affirmative answer is equivalent to the non-existence
of nontrivial evolutions of reduced algebras in equi-characteristic zero.

\vskip .1cm
\noindent{\bf Problem.\ }  Suppose $f\in {\bf C}[[x_1,\dots,x_n]]$
is a power series without constant term over the complex
numbers, and $I$ is the ideal of the reduced singular locus of
$f$, that is, $I$ is
the radical of the ideal generated by the partial derivatives
of $f$.  Does it follow that $f\in (x_1,\dots,x_n)I$?
\vskip .1cm

It is evident that the
answer is ``yes'' if $f$ has an isolated singularity.  Thus the
first interesting case occurs when the singular locus
of $f$ is a curve in ${\bf C}^3$.  The results of this paper
suffice to prove that the answer is still ``yes'' in this case,
or more generally whenever the embedding dimension of the reduced
singular locus of $f$ is less than $4$, or in the case of embedding
dimension $4$ when the reduced singular locus of f is Gorenstein, or
licci. In particular, the answer is ``yes'' for power series
in at most 3 variables.  We do not know any example where the
answer is ``no'' in characteristic 0 (or for the analogous
question in mixed characteristic).

In characteristic 0, when $f$ is quasihomogeneous,
one can express $f$ itself as a linear combination of the
derivatives of $f$, with coefficients
in $(x_1,\dots,x_n)$; thus in characteristic 0 the answer to
the question above is ``yes'' for any quasihomogeneous
polynomial.
In contrast, the answer to the question above
is ``no'' for the following  quasihomogenous polynomial in four variables over
a field of characteristic $p$, from the example treated at the
end of this paper:
$$
f= x_1^{p+1}x_2-x_2^{p+1}-x_1x_3^p+x_4^p.
$$

When $T$ is a complete
$\Lambda$-algebra that is a homomorphic image of
a power series ring of the form
$\Lambda[[x_1,...,x_n]]$, we could define evolutions similarly, using the
module of
universally finite derivations.  We then could extend the
notion to the non-local case by requiring the
condition above after localization or completion.
All the results of this paper could be cast in
that general setting.  But to simplify our presentation,
we shall work with local algebras that are essentially of finite type,
as indicated.

We began to study evolutions because of a situation that occurs in
the deformation theory of Galois representations.
In this theory one deals with a
``universal deformation ring" $R$, about which almost nothing
is known beyond the fact that it is a complete noetherian local ring, and
one also has a particular quotient of that ring, $T$,
a completion of a Hecke algebra, which is somewhat
more accessible.  From first order infinitesmal information about
deformations, which may be computed by cohomology, one
can sometimes show that the mapping
$R \to T$ is an evolution;
such arguments occur for example in Mazur [1994],
which is an exegesis of the work of Flach.
At this point if, by virtue of some
special properties of $T$, one knows that it admits no nontrivial
evolutions  (e.g., if T is a reduced complete intersection) then
one may deduce that $T$ is the universal deformation ring.
Such an assertian implies that certain Galois
representations are modular, and is interesting for that reason.

The work of Wiles and Taylor-Wiles uses other methods
to establish that the map  $R\to T$ is an isomorphism for the
rings $R$ and $T$ which they treat. Their methods show,
in addition, that their rings are reduced complete intersections.
This implies,
of course, that $T$ has no nontrivial evolutions, and gives
an a posteriori ``explanation'' why $R=T$.
We would like to have a clearer
understanding of evolutions
as a tool to be used in this type of argument, even in
cases where the rings in question are not complete intersections, or are
not known to be.

We are grateful to E. Kunz and H. Lenstra for sharing their work
with us and allowing us to include some of their results and
examples, to S. Popescu for his help, mathematical and
computer-related, with the examples,
and to C. Huneke,  and M. Hochster with whom we had
helpful discussions of the material.

\vskip .3cm

\qquad {\bf 1. Criteria for the existence of evolutions}
\vskip .2cm

We shall use a  criterion for the existence of evolutions
suggested by Lenstra.  If $T$ is a ring and $\phi: M \to N$ is an epimorphism
of
$T$-modules,
then $\phi$ is {\it minimal} if there is no proper submodule
$M'\subset M$
such that $\phi(M') = N$.

\proclaim Proposition \the\tno\ (Lenstra).
 Let $\Lambda$ be a Noetherian ring and
let $T$ be a local $\Lambda$-algebra, essentially
of finite type over $\Lambda$.
Every evolution
of $T$ is trivial
iff for some (equivalently all) presentations $T = P/I$, where
$P$ is a localization of a polynomial ring over $\Lambda$,
the map
$$
d_{T/P/\Lambda}: I/I^2 \to {\rm ker}(T\otimes_{P}\Omega_{P/\Lambda}
\to \Omega_{T/\Lambda})
$$
induced by the universal derivation is minimal.

\advance\tno by 1

\noindent {\sl Proof.\ } Let $T = P/I$ be any presentation of
$T$ where $P$ is a localization of a polynomial ring in finitely many
indeterminates over $\Lambda$.  If $J$ is an ideal
of $P$ with $J \subset I$ then an obvious diagram chase
shows that the natural surjection
$R := P/J \to P/I = T$ is an evolution iff
the differential
$
d_{T/P/\Lambda}
$
carries
$(J+I^2)/I^2$ onto the same image as $I/I^2$.
Using Nakayama's Lemma we see that $J = I$
iff $(J+I^2)/I^2 = I/I^2$, so $d_{T/P/\Lambda}$ is minimal iff
$T$ has no nontrivial evolution of the form $P/J$.

It remains to show that the condition that $d_{T/P/\Lambda}$ be
minimal is independent of the presentation $T = P/I$ chosen.
Since the family of presentations is filtered,
it is enough to show that if $T = P/I$ is a presentation,
$P'$ is a localization of a polynomial ring in one variable $x$ over $P$, and
$T = P'/I'$ is a presentation extending $T = P/I$ in an obvious sense,
then
$d_{T/P'/\Lambda}$ is minimal iff $d_{T/P/\Lambda}$ is minimal.

Let $g\in P$ be an element with the same
image in $T$ as $x$, so that $x-g\in I'$.  Replacing $x$ with the
new ``variable'' $x-g$, we may assume for simplicity that $g$ is 0.
We then have $I'/I'^2 = I/I^2\oplus Tx$ and
${\rm ker}(d_{T/P'/\Lambda})
= Tdx\oplus{\rm ker}(d_{T/P/\Lambda})$.
The required equivalence follows by Nakayama's Lemma.
\Box

Here is the relation of the problem given in the introduction
to the existence of evolutions in characteristic 0:

\vskip .2cm
\proclaim {Corollary} \the\tno.
There exists a reduced local ${\bf C}$-algebra $T$
of finite type whose localization at the origin has a nontrivial
evolution iff there exists a power series
$f\in {\bf C}[[x_1,\dots,x_n]]$ without constant term
such that
$$
f\notin (x_1,\dots,x_n)\sqrt{(f, df/dx_1,\dots,df/dx_n)}.
$$

\advance\tno by 1

\vskip .2cm

\noindent {\sl Proof.\ }  Set $I:= \sqrt{( f, df/dx_1,\dots,df/dx_n)}$
(since we are in characteristic 0 we could leave out f without
changing this ideal).  First suppose that
$
f\notin (x_1,\dots,x_n)I
$.
Replacing $f$ by an approximation
to very high order, we may suppose that $f$ is a polynomial.
Since
$
f\notin (x_1,\dots,x_n)I
$.
we may choose an ideal
$J\subset I$ generated
by polynomials
such that $J\neq I$ but $(J,f) = I$.
Writing P for the localization at the origin of the polynomial ring
in the $x_i$ it follows from the
definition that $P/J \to P/I$ is a nontrivial evolution.

Conversely, suppose $P/J \to P/I$ is a nontrivial evolution,
where $I$ is a radical ideal, and $P$ is the localization
of a polynomial ring over ${\bf C}$ at the origin.
Let $f'$ be a minimal generator
of $I$ that is not contained in $J$.  By definition,
$df'\in\Omega_{P/{\bf C}} $ goes to zero in
$(\Omega_{P/{\bf C}}\otimes P/I)/dJ$. Thus there is are
elements $g_i\in J$ and $p_i\in P$ such that
$d(f')-\sum r_id(g_i))$ is in
$I\Omega_{P/{\bf C}}$.  Let
$f := f'-\sum r_ig_i$, and note that $f$ is again a minimal
generator of $I$ which is not in $J$.
Since $g_i \in J\subset I$
we have
$$
df = df' - \sum r_i d(g_i) - \sum g_i d(r_i)
\in I\Omega_{P/{\bf C}}
$$
Thus the partial derivatives
of $f$ are all contained in $I$.  Since $I$ is radical
and $f = \notin {\cal M}I$ we obtain
$f\notin {\cal M}\sqrt{(f,df/dx_1,\dots,df/dx_n)}$ as required.
\Box

In special cases we can identify the kernel of the map
$
d_{T/P/\Lambda}
$
defined in Proposition 1.
If $I$ is an
ideal in a ring $P$ we define the $n$th {\it symbolic power} of $I$
to be the ideal
$$
I^{(n)} = \{f\in P| f\in I^n_Q \subset P_Q
{\rm\ for\ all\ minimal\ primes\ } Q {\rm\ of\ }I\}.
$$
Note that the symbolic powers depend only
on the isolated ($\equiv$ non-embedded) primary components of $I$, and
in particular the first symbolic power is the intersection of
the isolated primary components of $I$.

\edef \crit {\the\tno}
\proclaim Theorem \the\tno.
 Let $\Lambda$ be a Noetherian regular ring. Let
$P,{\cal M}$ be a localization of a
polynomial ring in finitely many variables
over $\Lambda$ and let $I$ be an ideal of $P$.
If $T:= P/I$ is reduced and generically separable over $\Lambda$,
then the kernel of
$
d: I/I^2 \to T\otimes_\Lambda \Omega_{P/\Lambda}
$
is $I^{(2)}/I^2$.
Thus every evolution of $T$ is trivial
iff
$I^{(2)}\subset {\cal M}I$.

\advance\tno by 1

\noindent {\sl Proof.\ } Let $Q$ be a minimal prime of $I$.
Since $\Lambda$ is regular
the prime $Q_{Q} \subset P_{Q}$ is a complete intersection.  By
hypothesis $P/Q$ is separable over $\Lambda$ and
thus $d: (I/I^2)_Q \to \Omega_{P_Q/\Lambda}$ is an injection.
Since $T\otimes_\Lambda \Omega_{P/\Lambda}$ is free over $T$
the map
$$
T\otimes_\Lambda \Omega_{P/\Lambda} \to
\oplus_Q T\otimes_{\Lambda}\Omega_{P_Q/\Lambda}
$$
is also a monomorphism, where the direct sum
is taken over the minimal primes
of $I$.  Thus the kernel of
$
d: I/I^2 \to T\otimes_\Lambda \Omega_{P/\Lambda}
$
is the same as the kernel of
$
I/I^2 \to
\oplus_Q (I/I^2)_Q,
$
which is $I^{(2)}/I^2$ as required.
\Box

This result shows that to prove that every evolution of a
reduced, generically separable algebra of finite type is
trivial it suffices to show that the symbolic square of the
defining ideal $I$ is contained in $I$ times the maximal ideal.
The rest of this paper is concerned with results of this type.

Before turning to these results we
mention a different approach, suggested by Huneke.
Scheja and Storch [1970] and B\"oger [1971]
prove under some additional hypotheses that the kernel of the
map
$
d: I/I^2 \to T\otimes_\Lambda \Omega_{P/\Lambda}
$
is $I'/I^2$, where $I'$ is the integral closure of $I^2$.  Rather
than studying symbolic squares, one might
study integral closures, and try to prove that
for some class of ideals $I$ we have $(I^2)' \subset {\cal M}I$.
The disadvantage of this approach is that in general the integral
closure is much smaller than the symbolic square, so presumably
$(I^2)' \subset {\cal M}I$ holds less
frequently than
$I^{(2)} \subset {\cal M}I$.  The potential advantage
is that questions of integral
closure can often be reduced to questions about zero-dimensional
ideals, which is generally not the case for questions about
the symbolic square.

\vskip .3cm

{\bf 2.  Symbolic powers and Fitting ideals}
\vskip .2cm

In this section we establish a simple connection
between Fitting invariants
and symbolic powers, and use it to analyze
symbolic powers in certain special cases.  First
some definitions:

We shall say that $I$ is
{\it unmixed\ } if all the associated primes of $I$ are isolated.
(Caution: The word unmixed
is sometimes used to denote an ideal all of
whose associated primes have the same dimension --- a stronger
condition.)

Fitting invariants are defined as follows:  Suppose
$M$ is a finitely generated
module over a commutative ring $R$, with a presentation
of the form $M \cong R^n/K$.  For each integer $i\geq 0$ the $i$th
Fitting ideal of M is defined to be
$F_i(M) = I_{n-i}(K)$, the ideal generated by the
$(n-1)\times (n-i)$ minors that can be formed from the elements
of $K$, regarded as vectors of length $n$.  The Fitting invariants
are independent of the presentation chosen, and commute with
base change.

We recall two easy facts about Fitting invariants:
Let $M$ be a finitely generated $R$-module.
If $M'\subset M$ is a submodule, then
$F_i(M) \subset F_i(M/M')$.  This is because given a presentation
of $M$ we may choose a presentation of $M/M'$ with the same
generators and more relations.  Also,
if $I\subset R$ is an ideal then
$F_i(M/IM) \subset F_i(M)+I$.
This is because we can derive a presentation for $M/IM$ from
one for $M$ by adding only the relations saying that elements
of $I$ times the generators are 0.

\vskip .2cm
\edef \thm2 {\the\tno}
\proclaim {Theorem} \the\tno.
If $I$ is an unmixed ideal of depth $\geq c$ in a Noetherian
ring and $x \in I^{(2)}$ then
$$
F_{c-1}(I/(x)) \subset I,
$$
If
$I$ is generically a complete intersection of codimension $c$
then the converse holds as well.

\vskip .2cm

\advance\tno by 1

\noindent {\sl Proof.\ }  First, suppose
that $x \in I^{(2)}$.  Since $I$ is unmixed, it suffices to
prove that
$
F_{c-1}(I/(x)) \subset I
$
after localizing at a minimal prime of $I$,
so we may assume that $x \in I^{2}$.  Thus $I/(x)$ surjects onto
$I/I^{2} = I\otimes R/I$, so
$
F_{c-1}(I/(x)) \subset F_{c-1}(I/I^{2}) \subset F_{c-1}(I)+I
$
and it suffices to show that  $F_{c-1}(I)\subset I$.  This
follows from the ``structure theorem''
of Buchsbaum-Eisenbud [1974].

Suppose that $I$ is generically a complete intersection
of depth $c$ and
$
F_{c-1}(I/(x))
$
is contained in $I$.
Since the associated primes of $I^{(2)}$ are all minimal
primes of $I$, we may begin by localizing
at one such and suppose that $I = (f_1,...,f_c)$ is a complete
intersection. Under these
circumstances we have $I^{2} = I^{(2)}$, and we shall prove
that $x\in I^{2}$.

The ideal $I$ is generated by the $c$ elements
 $f_i$.  We may take the same generators
for $I/(x)$. One of the relations on $I/(x)$ may be represented
as a column vector whose entries $g_i$ satisfy
$x = \sum g_i f_i$.
{}From the definition we see that $F_{c-1}(I/(x))$ contains
each $g_i$
so  we have $g_i\in I$ and
$x =  \sum g_i f_i \in I^{2}$.
\Box
\vskip .3cm

\noindent{\bf Remark.\ } In the case where
$I$ is generically a complete intersection,
the idea used in the
second part gives the following more general version:

\vskip .2cm
\proclaim {Theorem $4'$}.
If $I$ is an unmixed ideal of depth $\geq c$ in a Noetherian
ring and $I$ is generically a complete
intersection then $x \in I^{(d+1)}$  iff
$$
F_{N-1}(I^{(d)}/(x)) \subset I,
$$
where $N = {c+d-1\choose d}$.
\Box

\vskip .3cm


As an application we can show that the symbolic
square of a grade 2 perfect ideal in a ring $R$
 (that is,
an ideal that contains
an $R$-sequence of length 2
and that has projective dimension 1 as an $R$-module)
cannot be too big if the ideal is generically a complete intersection.
By the Hilbert-Burch theorem any such ideal can be represented
as the minors of an $n\times (n-1)$ matrix.

\vskip .2cm
\edef \theorem4 {\the\tno}
\proclaim {Corollary} \the\tno.
Suppose that an ideal $I$
in a ring $R$ is the ideal of $(n-1)\times (n-1)$ minors of an
$n\times (n-1)$ matrix $M$, that the depth of $I$ on
$R$ is two, and that $I$ is generically
a complete intersection in $R$.
If $J$ is the ideal generated by the entries of
any column of the matrix $M$, then $I^{(2)} \subset IJ$.

\vskip .2cm

\advance\tno by 1

\noindent {\sl Proof.\ } By the Hilbert-Burch Theorem,
a free resolution of $I$
can be written in the form
$$
0 \to R^{n-1}\to R^n \to I\to 0
$$
where the left hand map is $M$.  The generator of $I$ coming
from the $j$th generator of $R^n$ is the minor of $M$ involving
all but the $j$th row. Considering the expansion
of a determinant along a column, we see that every
element  $x\in I$ can be written as the determinant of
an $n\times n$ matrix $N$ whose first $n-1$
columns are the columns of $M$.  The matrix $N$ is also the
presentation matrix of $I/(x)$.  By
Theorem \thm2 ,
$x \in I^{(2)}$
iff $F_1(I/(x)) \subset I$, that is, iff the $(n-1)\times
(n-1)$ minors of $N$ are contained in $I$.  Expanding the
determinant of $N$ along a column of $M$ regarded as a column
of $N$, we see that $x$, the determinant, is in $IJ$ as
claimed.
\Box

\vskip .2cm
In the next section we
shall use a theorem of Buchweitz to extend
this result to any licci ideal that is generically a
complete intersection.
\vskip .3cm

Suppose that $I$ is an ideal
of depth $c$ in a ring $R$, and let $r \in F_c(I)$.
In any localization of $R[r^{-1}]$ at a prime ideal we can
generate $I$ by $c$ elements, so $I' := IR[r^{-1}]$
is locally a complete intersection, and
$I'^{(d)} = I'^d$ for all $d$.  This shows that a power of $r$,
and thus also a power of $F_c(I)$, annihilates $I^{(d)}/I^d$
for all $d$.  Using Theorem \thm2 \ we may sharpen this result for
$d=2$:

\vskip .2cm
\edef \Theorem3 {\the\tno}
\proclaim {Theorem} \the\tno.
Let $R$ be a Noetherian ring.  If
$I \subset R$
is an ideal that
is generically a complete intersection of codimension $c$,
then $F_c(I)$ annihilates $I^{(2)}/I^2$.

\advance\tno by 1

\vskip .2cm

\noindent {\sl Proof.\ }
Let $f_1,...,f_n$ generate $I$, and let $M:R^m\to R^n$ be
the matrix of relations on the $f_i$, that is, the
image of $M$ is the kernel of the map $f: R^n\to R$
defined by the $f_i$. Let $\gamma = (g_1,\ldots,g_n) \in R^n$, and write
$g = f(\gamma) = \sum g_if_i$ for the image of $\gamma$ in $I$.
Theorem \thm2 \
shows that $g$ is in
$I^{(2)}$ iff the matrix $M'$ obtained from $M$ by
adjoining a column with entries $g_i$ has all its
$(n-c)\times (n-c)$ minors in $I$. Since the
$(n-c)\times (n-c)$ minors of $M$ are already in $I$ by
the theorem of Buchsbaum-Eisenbud [1974],
this condition may be expressed by saying that the
minors of $M'$ containing the column of $g$'s are in $I$.

We may express this condition invariantly as follows.
The map $\wedge^{n-c-1}M$ corresponds to an element
$$
M_{n-c-1}\in \wedge^{n-c-1}R^n\otimes\wedge^{n-c-1}R^{m*},
$$
which can be written in terms of the minors of order $n-c-1$ of $M$.
The wedge product with this element defines a
map
$$
m: R^n\longrightarrow \wedge^{n-c}R^n\otimes\wedge^{n-c-1}R^{m*}.
$$
This map takes the element $\gamma$ to the vector whose entries
are the minors of $M'$ of order $n-c$ that involve the column
of $g_i$. Thus $g = f(\gamma) \in I^{(2)}$ iff
$$
m(\gamma) \in I\wedge^{n-c}R^n\otimes\wedge^{n-c-1}R^{m*}.
$$
Equivalently, the kernel of $R/I \otimes m$ maps by $f$ onto
$I^{(2)}/I$.

If
$\gamma$ is in the image of $M$ then $f(\gamma) = 0 \in I^{(2)}$,
so $m(\gamma)$ has entries in $I$.
Thus we obtain a complex
$$
(R/I)^m
\longrightarrow
(R/I)^n
\longrightarrow
\wedge^{n-c}(R/I)^n\otimes\wedge^{n-c-1}(R/I)^{m*}.
$$
where the left hand map is $R/I \otimes M$.  The cokernel
of $R/I \otimes M$ is $I/I^2$.  Thus the discussion
above may be summarized by saying that the
homology of this complex is $I^{(2)}/I^2$.

We shall conclude the proof by showing
that if $M$ is any $n\times m$ matrix over a ring $S$
such that the $(n-c)\times (n-c)$ minors of $M$ are zero,
then the complex derived from $M$ as in the display above,
replacing $R/I$ by $S$, has
homology annihilated by the ideal $J$ of
$(n-c-1)\times (n-c-1)$ minors of $M$; in fact
we shall show that multiplication by any such minor
is homotopic to 0 on this complex.  Let $M'$ be an
$n\times (m-c-1)$ submatrix of $M$; it suffices to show that
multiplication by
any maximal minor of $M'$ is homotopic to 0 on the subcomplex
$$
S^{m-c-1}
\longrightarrow
S^n
\longrightarrow
\wedge^{n-c}S^n\otimes\wedge^{n-c-1}S^{m-c-1*}.
$$
The dual of this complex is the right hand end of the
complex of Buchsbaum-Rim [1964].  In the case where
$M'$ is a matrix
whose entries are distinct indeterminates, the Buchsbaum-Rim
complex is a free resolution of the cokernel of
$M'^*$. The ideal $J$ then annihilates the cokernel of $M'$
(and is in fact exactly the annihilator of $coker(M'^*)$).
Thus multiplication by an $(n-c)\times (n-c)$
minor is homotopic
to 0 on the Buchsbaum-Rim complex in this case, and by
specialization this holds in every case.  (Buchsbaum and Rim
actually construct homotopies for these minors explicitly).
\Box

\vskip .3cm

It would be interesting
to know more about the annihilators of
$I^{(d)}/I^d$.  One simple case might be that
where $I$ is an unmixed codimension 2 homogeneous ideal in
a polynomial ring in 3 variables.
A preliminary study of examples produced by the program
Macaulay of Bayer and Stillman [1981]
suggests the following:

\vskip .2cm
\noindent{\bf Conjecture.\ } If
$I$ is the ideal of minors of a
``generic'' (that is, random) $2\times 3$
matrix of linear forms in 3 variables, then
the annihilator of
$I^{(d)}/I^d$
is
$
F_1(I)^e,
$
where $e$ is the greatest integer $\leq d/2$.

\vskip .5cm

\noindent{\bf 3.  Special cases}
\vskip .3cm
\quad\quad{\bf A. Monomial ideals}
\vskip .2cm

\edef \Theorem31 {\the\tno}
\proclaim {Proposition} \the\tno.
Suppose that $I$ is a monomial ideal in a polynomial ring
$k[x_1,\dots,x_n]$.  If $P$ is a monomial prime ideal
containing $I$ then for any $d\geq 0$ we have
$I^{(d)}\subset PI^{(d-1)}$.
\vskip .2cm

\advance\tno by 1

\noindent {\sl Proof.\ } Suppose that $I=\cap Q_j$ is a primary
decomposition.  A monomial ideal $Q_i$ is
primary, say to the monomial prime
$P_i = (x_{i_1},\dots x_{i_s})$,
iff $Q_i$ contains a power of each of the variables
$x_{i_t}$ and the minimal generators of $Q_i$ do not involve
any variables other than
$x_{i_1},\dots x_{i_s}$.
Thus any power of a primary monomial ideal is again primary,
and we have
$I^{(d)} = \cap Q_j^d$.

Now suppose that a monomial $m$ is in $I^{(d)}$.  By the argument above,
for each index
$j$ we may write $m = r_jm_{j,1}\cdots m_{j,d}$ with
$m_{j,i} \in Q_j$. Since $m\in P$, some variable $x_t \in P$
divides $m$, and thus divides (at least) one of
$r_j,m_{j,1},\dots, m_{j,d}$.  It follows that $m/x_t$ may be
written as a product with at least $d-1$ factors in $Q_j$,
so $m/x_t \in \cap Q_j^{d-1} = I^{(d-1)}$.  Thus $m \in PI^{(d-1)}$
as claimed.
\Box

\vskip .3cm
\quad\quad{\bf B. Licci ideals}
\vskip .2cm

Recall that a perfect ideal $I$ in a regular ring $R$ is called
{\bf licci\ } if it is in the linkage class of a complete intersection;
see for example Huneke-Ulrich [1987,1992] for some of the marvelous
properties of these ideals.  Licci ideals include for example all
perfect ideals of codimension 2 and all Gorenstein ideals of codimension
3.  Thus the following result greatly generalizes Corollary \theorem4 .

\vskip .2cm
\edef \Theorem32 {\the\tno}
\proclaim {Theorem} \the\tno.
Suppose that $R$ is a regular local ring, and that $I\subset R$ is
a perfect ideal that is generically a complete intersection.  Let
$J$ be the ideal generated by the elements in a row of some
presentation matrix (over $R$ or over $R/I$) of the canonical module
$\omega_{R/I}$ of $R$.  If $I$ is licci, then $I^{(2)}\subset JI$.
\vskip .2cm

\advance\tno by 1
\noindent {\sl Proof.\ }  An ideal $J$ as in the Theorem
is just the annihilator of some cyclic
quotient module $R/J \cong \omega_{R/I}/\omega'$
for some submodule $\omega' $ of $\omega$.  In particular,
such an ideal $J$ must contain $I$.
By a theorem of Buchweitz (Thesis, [1981]; see also
Buchweitz-Ulrich [1983]), the fact that $I$ is licci implies that
the module $\omega_{R/I}\otimes I/I^2$ is a Cohen-Macaulay module
over $R/I$.  Thus the map
$$
\omega_{R/I}\otimes I^{(2)}/I^2 \to \omega_{R/I}\otimes I/I^2
$$
is zero.  The natural map
$I^{(2)} \to I \to I/IJ$ can be factored as
$$
I^{(2)} \to R/J\otimes I^{(2)}/I^2 \to R/J\otimes I/I^2 = I/IJ
$$
and is thus 0, whence $I^{(2)}\subset IJ$.
\Box

\vskip .3cm
\quad\quad{\bf C. Almost complete intersections (Kunz)}
\vskip .2cm

Another case in which symbolic squares behave well is for
almost complete intersections.
The following result and its proof
was communicated to us by E.~Kunz (in a slightly
different form).  We
are grateful to him for allowing us to include it.

\vskip .2cm
\proclaim {Theorem (E.~Kunz)} \the\tno.
Let $P$  be a regular local ring with maximal ideal
be a ${\cal M}$, and let $I\subset P$
be an ideal that is generically a complete intersection.
If $I$ can be generated by $codim(I)+1$
elements, then $I^{(2)}\subseteq {\cal M}I$.

\advance\tno by 1
\vskip .2cm

\noindent {\sl Proof.\ }
Set $n := codim(I)$.  If
$I^{(2)}\not\subset{\cal M}I$, then $(I/I^{(2)})$
can be generated by $n$ elements.
Since $I/I^{(2)}$ is free of rank $n$ locally at any
minimal prime of $I$, we see that
$I/I^{(2)}\cong (P/I)^n$, and hence
$I/I^2\cong (P/I)^n\oplus N$ where $N$ is a nonzero module.
By the proposition of
Vasconcelos [1968] $I$ contains a
regular sequence of length greater than $n$,
a contradiction.
\Box

\vskip .3cm

\vskip .3cm
\quad\quad{\bf D. Quasi-homogeneous ideals}

\vskip .2cm
\edef \tag {\the\tno}
\proclaim {Proposition} \the\tno.
Let $P = k[x_1,\dots,x_n]$
be a polynomial ring over a field $k$, and let
${\cal M} = (x_1,\dots,x_n).$ Suppose that
$I\subset P$ is an unmixed ideal which is quasihomogeneous (that is,
homogeneous with respect to some system of strictly positive
integer weights of the variables).
If $f\in I^{(d)}$ is a quasihomogeneous
element, then $f\in deg(f)\cdot {\cal M}I^{(d-1)}$.  In particular,
if $char(k) = 0$ then
$I^{(d)}\subseteq {\cal M}I^{(d-1)}$,
and if in addition $I$ is a radical ideal then every evolution
of $P/I$ is trivial.

\advance\tno by 1
\vskip .2cm

\noindent {\sl Proof.\ }  If $f\in I^{(d)}$ then
by definition there is a
polynomial $h$, not contained in any of the minimal primes of $I$,
 such that
$hf \in I^d$.
Differentiating we see that for each index $j$ we have
$hdf/dx_j+fdh/dx_j \in I^{d-1}$.  Since $f\in I^{(d-1)}$
we deduce
$hdf/dx_j\in I^{(d-1)}$, whence
$df/dx_j\in I^{(d-1)}$.  Euler's relation
$deg(f)\cdot f = \sum deg(x_j)\cdot x_jdf/dx_j$ now shows that
$f\in deg(f)\cdot {\cal M}I^{(d-1)}$ as required.  The last assertion
follows from Proposition 1.
\Box

\vskip .3cm

In positive characteristic, by contrast, there are evolutions of
quasihomogeneous rings, even of 1-dimensional
quasihomogeneous domains.  The first example of which we are aware
was communicated to us by E.~Kunz.  It is the (localization at the
quasihomogeneous maximal ideal of the) domain
$$
k[t^{14},t^{20}, t^{25},t^{30}, t^{91}],
$$
where $k$ is a field of characteristic 2.
A careful analysis of Kunz' example led us
to a series of simpler examples, in all positive characteristics.
We would like to thank Sorin Popescu, who helped us
greatly, both with
the computation (using Maple and the program Macaulay of
Bayer and Stillman [1981--]) and the mathematical aspect
of the examples below.

These examples are in certain senses minimal:
It follows from the results above that any example of a reduced
algebra with nontrivial evolution must have
embedding codimension at least 3, and our examples have
embedding codimension exactly 3.  Also, any local Cohen-Macaulay
ring with
embedding codimension at least 3 has
multiplicity at least 4.  Our examples are one-dimensional
domains, thus Cohen-Macaulay,
and the example in characteristic 2 below
has multiplicity exactly 4.

\vskip .2cm
\noindent{\bf Example}
Let $k$ be a field of characteristic $p>0$, and let $I$ be the
kernel of the map
$$\eqalign{
k[x_1,\dots,x_4] &\to k[t]\cr
x_1,x_2,x_3,x_4 &\mapsto
t^{p^2},t^{p(p+1)}, t^{p^2+p+1},t^{(p+1)^2}}
$$
(or the localization of this ideal at the maximal ideal
$(x_1,\dots,x_4)$).
Let
$$
f= x_1^{p+1}x_2-x_2^{p+1}-x_1x_3^p+x_4^p
$$
We claim that $f$ is a minimal generator of $I$ but $f$ is
contained in $I^{(2)}$.  (It is easy to see that the
derivatives of $f$ are contained in $I$, but this is
not quite equivalent in characteristic $p$.)

To show that $f\in I^{(2)}$,
consider
the polynomials
$$\eqalign{
g_1 &= x_1^{p+1}-x_2^p     \cr
g_2 &= x_1x_4-x_2x_3    \cr
g_3 &= x_1^px_2-x_3^p.  \cr
}
$$
One checks immediately by applying the homomorphism
of rings above that $f, g_1,\dots,g_3 \in I$.  Since $I$
is prime, the relations
$
x_1^pf = g_1g_3+g_2^p
$
and
$
x_1 \notin I
$
show that $f\in I^{(2)}$.

To prove that $f$ is a minimal generator of $I$ it suffices
to show that no element of $I$ has a term of the form
$x_4^a$ with $0<a<p$.  Since $I$ is generated by binomials,
it suffices to show that there is no binomial of the form
$x_4^a-x_3^bx_2^cx_1^d$ in $I$, or equivalently that the
equation
$$
a(p+1)^2 = b(p^2+p+1) + cp(p+1) + dp^2
$$
cannot be satisfied by nonnegative integers $a,b,c,d$ with
$0<a<p$.

This is elementary:  Suppose $0< a < p$.
Reducing modulo $p$ we see that
$a\equiv b \ ({\rm mod}\ p)$.  If $b = a+np$ for some $n\geq 1$
then subtracting $a(p+1)^2$ from both sides of the
displayed relation we get
$
(p-1)p\geq ap \geq np(p^2+p+1),
$
which is impossible.
Thus $b=a$.  Subtracting
$
a(p^2+p+1)
$
from both sides of the
displayed relation
we get
$
p(p-1)\geq ap = cp(p+1)+dp^2
$
and the right hand side is either 0 or $\geq p^2$,
a contradiction.
\vskip .3cm

{
\baselineskip=12pt

\centerline {\bf References}

\vskip .1cm
\frenchspacing


\item{} D.~Bayer and  M.~Stillman.
Macaulay, a computer algebra system.  Available free from the authors.
(1981-).
\vskip .1cm

 \item{} E.~B\"oger.
Differentielle und ganz-algebraische
Abhangigkeit bei Idealen analytischer Algebren.
{\sl Math. Z.} {\bf 121} (1971)  188--189.
\vskip .1cm

\item{} D.~A.~Buchsbaum and D.~Eisenbud.
Some structure theorems for finite free resolutions.
{\sl Adv.~in Math.~} {\bf 12} (1974), 84--139.
\vskip .1cm

\item{} D.~A.~Buchsbaum and D.~S.~Rim,
A generalized Koszul complex. II. Depth and multiplicity.
{\sl Trans. Amer. Math. Soc.} {\bf 111} (1964)  197--224.
\vskip .1cm

\item{} R.-O.~Buchweitz: Contributions \`a la th\'eorie
des singularit\'es:
1. D\'eformations de Diagrammes, D\'eploiements et
        Singularit\'es tr\`es rigides, Liaison alg\'ebriques.

2. On Zariski's Criterion for Equisingularity and Non-Smoothable
        Monomial Curves

Th\`ese, Paris VII (1981).
\vskip .1cm

\item{} R.-O.~Buchweitz, B.~Ulrich:
Homological Properties which are invariant under Linkage.
preprint Purdue (1983).
\vskip .1cm

\item{} D.~Eisenbud and M.~Hochster.
 A Nullstellensatz with
nilpotents, and Zariski's Main Lemma on holomorphic functions.
{\sl J. Alg.} {\bf 58} (1979) 157--161.
\vskip .1cm

\item{} M.~Flach.
A finiteness theorem for the symmetric square of an elliptic curve.
{\sl Invent. Math.} {\bf 109} (1992) 307--327.
\vskip .1cm

\item{} J.~Herzog.
Ein Cohen-Macaulay Criterium mit
Anwendungen auf den Konormelen-modul und den Differentialmodul.
{\sl Math. Z.} {\bf 163} (1978) 149--162.
\vskip .1cm

\item{} C.~Huneke, and B.~Ulrich.
The structure of linkage.
{\sl Ann. of Math.} (2) {\bf 126} (1987) 277--334.
\vskip .1cm

\item{} C.~Huneke, and B.~Ulrich.
Local properties of licci ideals.
{\sl Math. Z.} {\bf 211} (1992) 129--154.
\vskip .1cm

\item{} A.~Kustin and B.~Ulrich.
A family of complexes associated to an almost alternating map, with
applications to residual intersections.
{\sl Mem. of the
Amer.~Math.~Soc.~}{\bf 461} (1992).
\vskip .1cm

\item{} B.~Mazur.
Deformations of Galois representations and Hecke Algebras.
Harvard Course Notes, available by request from the
author (1994).
\vskip .1cm

\item{} M.~Nagata.
{\bf Local Rings}.
Wiley Interscience, New York 1962.
\vskip .1cm

\item{} K.~Saito.
Quasihomogene isolierte Singularit\"aten von Hyperfl\"achen.
{\sl Invent.~Math.~} {\bf 14} (1971) 123--142.
\vskip .1cm

\item{} G.~Scheja and U.~Storch.
\"Uber differentielle
Abh\"angigkeit bei Idealen analytischer Algebren.
{\sl Math. Z.} {\bf 114} (1970) 101--112.
\vskip .1cm

\item{} W.~Vasconcelos.
A note on normality and module of differentials.
{\sl Math. Z.} {\bf 105} (1968)  291--293.
}
\vskip .4cm
\noindent{Dept. of Math., Brandeis Univ., Waltham MA 02254,
eisenbud@math.brandeis.edu}
\vskip .1 cm
\noindent{Dept. of Math., Harvard Univ., Cambridge MA 20138,
mazur@math.harvard.edu}

\bye